\documentclass{svmult}
\usepackage[hidelinks]{hyperref}
\usepackage{mathtools} 
\usepackage{type1cm}        
\usepackage{makeidx}         
\usepackage{graphicx}        
\usepackage{multicol}        
\usepackage[bottom]{footmisc}

\usepackage{newtxtext}       %
\usepackage[varvw]{newtxmath}       

\makeindex             
\usepackage{bbm}

\usepackage{microtype}
\usepackage{subcaption}
\usepackage{tikz}
\usepackage{enumitem}
\usepackage{mathrsfs}
\usepackage[utf8]{inputenc}
\usepackage[style=alphabetic, backref=false, giveninits=true,citestyle=alphabetic, natbib=false, url=false, isbn=false,doi=true,eprint=false,maxbibnames=99]{biblatex}
\DeclareFieldFormat[article]{title}{{#1}}
\DeclareFieldFormat[unpublished]{title}{{#1}}
\DeclareFieldFormat[proceedings]{title}{{#1}}
\DeclareFieldFormat[incollection]{title}{{#1}}
\DeclareFieldFormat[unpublished]{title}{{\em #1}}
\renewbibmacro{in:}{}
\DeclareFieldFormat[article]{pages}{#1}
\spnewtheorem{Exa}[theorem]{Example}{\normalfont\bfseries}{\normalfont}
\spnewtheorem{assumptions}[theorem]{Assumptions}{\normalfont\bfseries}{\normalfont}
\spnewtheorem{Prop}[theorem]{Proposition}{\normalfont\bfseries}{\normalfont}
\spnewtheorem{Lem}[theorem]{Lemma}{\normalfont\bfseries}{\normalfont\itshape}
\spnewtheorem{Fact}[theorem]{Fact}{\normalfont\bfseries}{\normalfont\itshape}
\spnewtheorem{KLem}{Key Lemma}{\normalfont\bfseries}{\normalfont\itshape}
\spnewtheorem{Cor}[theorem]{Corollary}{\normalfont\bfseries}{\normalfont\itshape}

\newcommand{\supp}{\operatorname{supp}}
\renewcommand{\D}{\mathcal{D}}

\DeclareMathOperator{\dom}{dom}

\renewcommand{\tilde}{\widetilde}
\renewcommand{\hat}{\widehat}
\renewcommand{\emptyset}{\varnothing}
\renewcommand{\d}{\;\mathrm{d}}

\newcommand{\dm}{\;\mathrm{d}m}
\newcommand{\1}{\mathbbm{1}}
\newcommand{\R}{\mathbb{R}}

\renewcommand{\S}{\mathcal{S}}

\begin{document}
\title*{Generalised Kre\u{\i}n--Feller operators and gap diffusions via transformations of measure spaces}
\titlerunning{Generalised Kre\u{\i}n--Feller operators and gap diffusion}
\author{Marc\ Kesseb\"ohmer, Aljoscha Niemann, Tony Samuel, and Hendrik Weyer}
\authorrunning{M.~Kesseb\"ohmer, A.~Niemann, T.~Samuel, H.~Weyer}
\institute{M.\ Kesseb\"ohmer, A.\ Niemann and H.\ Weyer \at FB 3, Mathematik und Informatik, University of Bremen, 28359 Bremen, Germany \email{mhk@uni-bremen.de, niemann1@uni-bremen.de, hendrik.weyer@gmail.com}
\and Tony Samuel \at School of Mathematics, University of Birmingham, Edgbaston, Birmingham, B15 2TT, UK \email{t.samuel@bham.ac.uk}}

\maketitle

\abstract{We consider the generalised Kre\u{\i}n--Feller operator $\Delta_{\nu, \mu} $ with respect to compactly supported Borel probability measures $\mu$ and $\nu$ with the natural restrictions that $\mu$ is atomless, the  $\supp(\nu)\subseteq\supp(\mu)$ and the atoms of $\nu $ are embedded in the $\supp(\mu)$. We show that the solutions of the eigenvalue problem for $\Delta_{\nu, \mu} $ can be transferred to the corresponding problem for the classical Kre\u{\i}n--Feller operator $\Delta_{\nu \circ F_{\mu}^{-1}, \Lambda}$ with respect to the Lebesgue measure $\Lambda$ via an isometric isomorphism determined by the distribution function $F_\mu$ of $\mu$. In this way, we obtain a new characterisation of the upper spectral dimension and consolidate many known results on the spectral asymptotics of  Kre\u{\i}n--Feller operators. We also recover known properties of and connections to generalised gap diffusions associated to these operators.}


\begin{center}
\textit{Dedicated to the memory of Professor Bob Strichartz.}
\end{center}

\section{Introduction and motivation}

The classical Kre\u{\i}n--Feller differential operator $\Delta_{\nu, \Lambda}$, where $\nu$ denotes a non-atomic compactly supported Borel probability measure on $\mathbb{R}$ and where $\Lambda$ denotes the \mbox{one-dimensional} Lebesgue measure, was introduced in \cite{Fe57,KK68}; in the case $\nu=\Lambda$, this operator coincides with the classical second order weak derivative. The spectral properties of the classical Kre\u{\i}n--Feller operator have been heavily investigated by, for instance, Hong and Uno \cite{MR118891}, McKean and Ray \cite{mckean_jr.1962}, Kac \cite{Kac_1973}, Fujita \cite{Fu87}, K\"uchler \cite{MR574035}, Langer \cite{MR0314125}, and Kotani and Watanabe \cite{MR661628}. In \cite{MR1328700} the result of Fujita \cite{Fu87} was refined by a Dirichlet form approach and renewal theory. Volkmer \cite{Volkmer05,MR2218192} was able to transform the Kre\u{\i}n--Feller eigenvalue problem to a semi-definite Sturm-Liouville problem and obtained explicit estimates of the eigenvalues in the case when $\nu$ is the Cantor measure. In the case when $\nu$ has no atoms, it has been established that $\Delta_{\nu, \Lambda}$ is the infinitesimal generator of a gap diffusion (also known as skip-free diffusion, quasi-diffusion or generalised diffusion), see for example \cite{Burkhardt1983,MR2817342,MR0345224,MR3034785,MR574035,MR0314125,MR3005002}.

Here, we investigate generalised Kre\u{\i}n--Feller operators $\Delta_{\nu, \mu}$ for Borel probability measures $\nu$ and $\mu$ on the real line with bounded support under the natural assumptions that $\mu$ is atomless and that  $\nu$ is {\em embedded} in  $\mu$, meaning $\supp(\nu)\subseteq\supp(\mu)$ and the possible atoms of $\nu$ are disjoint from all boundary points of the individual open connected components of $\R\setminus\supp(\mu)$. The generalised Kre\u{\i}n--Feller operator, has been introduced by K\"uchler \cite{MR574035} in the case that the distribution functions of $\mu$ and $\nu$ are strictly increasing and continuous, for general atomless measures by Freiberg and Z\"ahle \cite{FZ02}, and for generalised signed measures (referred to as the \textsl{measure Sturm-Liouville problem}) by Volkmer \cite{Volkmer05}. The associated Dirichlet forms have also been studied in \cite{ChenFukushima+2011}. Spectral properties of generalised Kre\u{\i}n--Feller operator have been considered by, for instance Arzt \cite{A15b}, Ehnes \cite{Ehnes2019}, Freiberg \cite{Fr05}, Rastegaev \cite{MR3407794}, She\u{\i}pak \cite{MR3370516}, Vladimirov and She\u{\i}pak \cite{MR3185121}, and Kesseb\"ohmer and Niemann \cite{KN21,KN22A}. For results on spectral problems for the classical Kre\u{\i}n--Feller operator in higher dimensions we refer the reader to \cite{MR1298682,MR2261337, Ngai_2021,KN21b,KN2022}.

In this article, generalising ideas from \cite{KSW16}, we elaborate on the connections between the generalised and the classical one-dimensional Kre\u{\i}n--Feller operators by establishing a suitable isometric isomorphism determined by the distribution function of $\mu$. This isometric isomorphism first appeared in \cite[Chapter 2]{MR2563669} under the additional assumption that the distribution function of $\mu$ (called scaling function) is strictly increasing on the convex hull of the support of $\mu$. Here, we extend this to the more general framework, assuming that $\mu$ is atomless and $\nu$ is embedded in $\mu$, and apply this connection to the associated spectral and stochastic properties. In doing so, we summarise and consolidate the results on spectral asymptotics scattered in the literature, for instance, in \cite{MR0278126,,Fr05,MR2787628,FZ02,KSW16,mckean_jr.1962,MR1328700,Volkmer05}, and give a new perspective concerning gap diffusions as studied, for example, in \cite{MR647781,MR574035,kuechler1986,MR0314125,MR1127331}.

As a first application of this observation we show that the spectral properties of the generalised Kre\u{\i}n--Feller operators can be reduced to those of the classical ones. We determine a general upper bound  for the spectral dimension in terms of the upper Minkowski dimension of the support of the involved measures, and deduce the strong spectral asymptotics for self-similar measures. Moreover, we obtain a new proof for those pairs of measures $\mu$ and $\nu$ for which $\nu$ has an absolutely continuous component with respect to $\mu$. As a second application, we connect properties of the associated gap diffusion for generalised Kre\u{\i}n--Feller operators to that of the classical \mbox{Kre\u{\i}n--Feller} operators. This complements and partially resembles the general framework established in \cite{DynkinI_II,Ogura_1989}.

\section{Set-up and summary of results}
Before we state our main results, we need the following basic setups and standing hypotheses as stated next. 
\subsection{Setup and standing hypotheses}\label{sec:standing_assumptions}

Let $\mu$ and $\nu$ denote two Borel probability measures on $\mathbb{R}$ with bounded support such that $\mu$ is atomless and $\nu$ is embedded in $\mu$. Without loss of generality we also assume that $\supp(\mu) \subseteq [0,1]$. Denote the distribution function of $\mu$ and $\nu$ by $F_{\mu}$ and $F_{\nu}$, respectively. Let $(C_{\nu, \mu}, \lVert \, \cdot \, \rVert_{\infty} )$ denote the Banach space of continuous functions with domain $[0,1]$ and which are linear in scale $F_{\mu}$ on intervals where $F_{\nu}$ is constant. Namely, on each connected component $J$ of $[0,1] \setminus \supp(\nu)$ the function $f$ is linear in scale $F_{\mu}$, that is $f(x)=a_J F_{\mu}(x)+b_J$ for all $x \in J$ and some $a_J$ and $b_J \in \mathbb{R}$. As indicated above, we let $\Lambda$ denote the one-dimensional Lebesgue measure restricted to $[0,1]$.

Set $\S^w \coloneqq L^2(\nu)$ and $\S^s \coloneqq C_{\nu,\mu}$, where $w$ stands for \textsl{weak} and $s$ stands for \textsl{strong}; for $* \in \{s,w\}$, we sometimes write $\S^*(\mu,\nu)$ instead of $\S^*$ to stress the dependence of the underlying measure spaces. In what follows, we will mainly be concerned with the Banach spaces $(\S^{w},\lVert\, \cdot \,\rVert_{L^2(\nu)})$ and $(\S^{s},\lVert\, \cdot \,\rVert_{\infty})$. Letting $* \in \{s,w\}$ be fixed, a function $f$ belonging to the set $C([0,1])$, of continuous functions with domain $[0,1]$, is said to lie in $\D^*(\Delta_{\nu,\mu})$ if there exist $a, b \in \mathbb{R}$ and $g \in \S^*$ with 
	\begin{align}\label{eq:KreinFeller}
	f(x)=a+bF_{\mu}(x)+\int_{[0,x]}( F_{\mu}(x)-F_{\mu}(y) ) g(y) \d \nu(y)
	\end{align}
for all $x \in [0,1]$. By Fubini's Theorem we have 
	\begin{align*}	
	f(x)=a+bF_{\mu}(x)+\int_{[0,x]}\int_{[0,y]} g(s)\d \nu(s)\d \mu(y).
	\end{align*}
If $F_{\mu}(x) \neq 0$, then by \eqref{eq:KreinFeller}
	\begin{align*}
	  \frac{f(x)-f(0)}{F_{\mu}(x)}&=b+\int_{[0,x]}\dfrac{F_{\mu}(x)-F_{\mu}(y)}{F_{\mu}(x)}g(y) \d \nu(y).
	\end{align*} 
	Since $\nu(\{0^+\})=0$ where $ 0^+ \coloneqq \inf (\supp(\mu))$, it follows 
		\begin{align*}\nabla_\mu f(0)\coloneqq
	 \lim_{x \searrow 0^+} \frac{f(x)-f(0)}{F_{\mu}(x)}&=b+\lim_{x \searrow 0^+} \int_{[0,x]}\dfrac{F_{\mu}(x)-F_{\mu}(y)}{F_{\mu}(x)}g(y) \d \nu(y)=b,
	\end{align*} 
	where we used that \begin{align*}
	  \int_{[0,x]}\dfrac{F_{\mu}(x)-F_{\mu}(y)}{F_{\mu}(x)}|g(y)|\d \nu(y) \leq 
	   \int_{[0,x]}|g(y)| \d \nu(y).
	\end{align*}
which converges to zero as $x$ approaches $0^{+}$. In particular, $b=\nabla_\mu f(0)$ depends only on $\mu$ and $f$. To see that $g$ is also uniquely determined by $f$, we assume the contrary. Namely, letting $* \in \{s,w\}$, we assume that there exist distinct $g_{1}$ and $g_{2} \in \S^*$ satisfying \eqref{eq:KreinFeller}. In which case, for $i \in \{1, 2\}$, we set
	\begin{align*}
	G_{i}(x) \coloneqq b+\int_{[0,x]} g_{i} \d \nu(y),
	\end{align*}
for all $x \in [0,1]$, and observe by Fubini's theorem 
	\begin{align}\label{eq:FTC}
	\begin{aligned}
	\int_{[0,x]} G_{i}(y) \d \mu(y)&=F_{\mu}(x) \nabla_{\mu}f(0)+
	\int_{[0,x]} \int_{[0,y]} g_{i}(z)\d \nu(z)\d\mu(y)\\
	&=F_{\mu}(x) \nabla_{\mu}f(0)+
	\int_{[0,x]}\int_{[z,x]} g_{i}(z) \d \mu(y)\d \nu(z)\\
	&=F_{\mu}(x) \nabla_{\mu}f(0)+
	\int_{[0,x]}(F_{\mu}(x)-F_{\mu}(z)) g_{i}(z)\d \nu(z)\\
	&=f(x)-f(0).
	\end{aligned}
	\end{align}
Using the condition that $\nu$ is embedded in $\mu$, which implies $G_i$ is continuous at the boundary points of all complementary intervals of $\supp(\mu)$, the right continuity of $G_{i}$, and the fact that $f$ is determined pointwise, one obtains a contradiction to the uniqueness of densities. This allows us to define the {\em (generalised) Kre\u{\i}n--Feller operator} $\Delta_{\nu,\mu}f \coloneqq g$ via \eqref{eq:KreinFeller}; we distinguish between the {\em strong} and {\em weak} operator depending on the chosen domain. Additionally, from this and \eqref{eq:KreinFeller}, one may conclude
\begin{align*}
\D^s ( \Delta_{\nu,\mu} ) \subseteq \D^w (\Delta_{\nu,\mu}) \subseteq C_{\nu,\mu}.
\end{align*}
The following example shows that the assumption, $\nu$ is embedded in $\mu$, cannot be removed in order to define $\Delta_{\nu,\mu}$. 
\begin{Exa}
If $\mu = 2\Lambda|_{[0,1/2]}$ and  $\nu$ is the Dirac measure with point mass at $1/2$, then \mbox{$\supp(\nu)\subseteq\supp(\mu)$,} but $\nu$ is not embedded in $\mu$. An elementary computation shows that the representation in \eqref{eq:KreinFeller} is independent of $g$; consequently $\Delta_{\nu,\mu}$ is not necessarily well defined for such measures.
\end{Exa}

For $\gamma =(\alpha,\beta) \in [0,{\pi/2}]^{2}$ we consider the eigenvalue problem of Fujita \cite{Fu87} with \textsl{Robin boundary conditions}; namely, to classify those $f \in \D^*(\Delta_{\nu,\mu})$ and $\lambda \in \mathbb{R}$ for which $\Delta_{\nu,\mu}f= \lambda f$ and
	\begin{align}\label{eq:EWC}
	 f(0)\cos(\alpha)-\nabla_{\mu}f(0)\sin(\alpha)=0,
	\;\;  f(1)\cos(\beta)+\nabla_{\mu}f(1)\sin(\beta)=0.
	\end{align}
The particular case $\gamma=(\pi/2,\pi/2)$ is referred to as the \textsl{Neumann case} and the case $\gamma=(0,0)$ is called the \textsl{Dirichlet case}. We denote by $\D^{*}_{\gamma}( \Delta_{\nu,\mu} )$ the set of $f \in \D^{*}( \Delta_{\nu,\mu} )$ which satisfy \eqref{eq:EWC}. Combining \eqref{eq:KreinFeller}~--~\eqref{eq:FTC} with Fubini's Theorem and our assumptions on $\nu$ and $\mu$ one obtains a \textsl{Gauss-Green formula}; namely, for $f$ and $g \in \D_{\gamma}^s( \Delta_{\nu,\mu} )$,
	\begin{align}\label{eq:gauss_green}
	\begin{aligned}
	&\int (\Delta_{\nu,\mu}f )g\d\nu\\
	&= (\nabla_{\mu}f (1)-\nabla_{\mu}f (0))g(0)+ \int\nabla_{\mu}g(y)(\nabla_{\mu}f (1)-\nabla_{\mu}f (y) ) \d\mu(y)\\
	&= \nabla_{\mu}f (1) g(1) - \nabla_{\mu}f (0) g(0) - \int \nabla_{\mu}f\,\nabla_{\mu} g \d\mu.
	\end{aligned}	
	\end{align}
 This in tandem with our boundary conditions implies that
 	\begin{align*}
	&\int (\Delta_{\nu,\mu}f )g \d \nu - \int (\Delta_{\nu,\mu}g )f \d\nu\\
	&\;\;\;\;\;= \nabla_{\mu}f (1) g(1) - \nabla_{\mu}f (0) g(0) - \nabla_{\mu}g(1) f(1) + \nabla_{\mu}g(0)f(0)= 0.
	\end{align*}
Hence, the operator $\Delta_{\nu,\mu}$ restricted to $\D^{*}_{\gamma}( \Delta_{\nu,\mu})$ is symmetric.  Further, setting $g = f$ in \eqref{eq:gauss_green}, and applying our boundary conditions once more, yields that $\Delta_{\mu,\nu}$ restricted to $\D^{*}_{\gamma}( \Delta_{\nu,\mu})$ is non-positive.

\subsection{Summary of results}

In Theorem \ref{thm:laplace_backward} we establish a strong connection between $\Delta_{\nu,\mu}$ and $\Delta_{\nu \circ F_{\mu}^{-1}, \Lambda}$. Indeed, by utilising the \textsl{pseudo-inverse} 
	\begin{align*}
	\check{F}^{-1}_{\mu}(x) \colon x \mapsto \inf\{ y \in [0,1] \colon F_{\mu}(y) \geq x \}
	\end{align*}
of $F_{\mu}$, we prove, for $* \in \{s,w\}$ and $\gamma \in [0,\pi/2]^{2}$, that $\varphi \colon f \mapsto f \circ \check{F}_{\mu}^{-1}$ is an isometric isomorphism on $\S^{*}$ with
	\begin{align*}
	\Delta_{\nu \circ F_{\mu}^{-1}, \Lambda }\circ\varphi=\varphi\circ \Delta_{\nu,\mu }\quad \text{and}\quad\varphi ( \D^*_{\gamma}(\Delta_{\nu,\mu}))=\D^*_{\gamma}(\Delta_{\nu \circ F_{\mu}^{-1},\Lambda}).
	\end{align*}
With this at hand, we are able to show the following which consolidates and extends known results.  In particular, that the spectral properties of $\Delta_{\nu \circ F_{\mu}^{-1},\Lambda}$ are inherited from $\Delta_{\nu ,\mu}$ and vice versa.

\begin{enumerate}
%
\item Theorem \ref{thm:USD} concerns the exponent of the asymptotic growth rate of the eigenvalue counting function of $\Delta_{\nu,\mu}$, namely the upper spectral dimension of $\Delta_{\nu,\mu}$. For a certain class of self-similar measures, in Theorem \ref{thm:Freiberg} we recover the asymptotic growth rate of the eigenvalue counting function of $\Delta_{\nu,\mu}$, which was first observed in \cite{Fr05}. Both of these results are achieved via an application of Theorem \ref{thm:laplace_backward}. 
Here, we also make use of \cite{KN22A} in deriving  the upper spectral dimension, and of \cite{MR1328700} to treat the spectral asymptotic growth rate in the self-similar setting.

\item Theorem \ref{thm:RayMacKean} concerns the asymptotic growth rate of the eigenvalue counting function of $\Delta_{\sigma^2 \mu+ \eta,\mu}$ where $\eta$ is an atomless Borel measure singular to $\mu$ and $\sigma^2  $ is a $\mu$-integrable function. This result is proved by combining Theorem \ref{thm:laplace_backward} and the fact that $\eta = \nu \circ F_{\mu}^{-1}$ is singular to the Lebesgue measure with the corresponding result for the case $\mu=\Lambda$ as given in \cite{MR0278126,mckean_jr.1962,Volkmer05}. Note, some of the proofs given in \cite{mckean_jr.1962} were pointed out to contain gaps by Kac \cite{Kac_1973}.
\item Letting $(X_{t})_{t \geq 0}$ denote the gap diffusion with speed measure $\nu \circ F_{\mu}^{-1}$,
utilising our correspondence theorem (Theorem \ref{thm:laplace_backward}), we show that $(\check{F}^{-1}_{\mu}( X_t ))_{ t\geq 0}$ is a Feller process (with respect to an appropriate topology on the state space) and its infinitesimal generator coincides with the generalised Kre\u{\i}n--Feller operator $\Delta_{\nu,\mu}$ with Neumann boundary condition;  
this gives a new outlook on \cite{MR647781,MR574035,kuechler1986,MR0314125,MR1127331,Ogura_1989}.

\end{enumerate}

\section{Kre\u{\i}n--Feller operators}

In this section, we collect some important properties  Kre\u{\i}n-Feller operators, starting with the classical one, where the reference measure is the Lebesgue measure. 
\subsection{Properties of classical Kre\u{\i}n--Feller operators} 

We now consider the case $\mu= \Lambda$, with respect to weak and strong solutions. Most of these results are nowadays folklore and can be found, for instance, in \cite{MR0314125}. Since we could not locate references where all the facts are proved in detail, here we give an overview and essentially reduce all properties down to two key observations, namely that, for $* \in \{ s, w \}$, the operator $\Delta_{\nu, \Lambda}$ restricted to $\D^{*}_{\gamma}( \Delta_{\nu,\Lambda})$ is symmetric (a consequence of the Gauss-Green formula), and that $\Delta_{\nu, \Lambda} \colon \D^{*}_{\gamma}( \Delta_{\nu,\Lambda}) \to \S_{\gamma}^*$ is surjective, as is shown in the following lemma. 

\begin{Lem}\label{lem:Surjec}
For $\gamma=(\alpha, \beta) \in [0,{\pi/2}] ^{2}$ and $* \in \{s,w\}$, the map 
\begin{align*}
\Delta_{\nu,\Lambda} \colon \D^{*}_{\gamma}( \Delta_{\nu,\Lambda} ) \rightarrow \S_{\gamma}^*
\end{align*}
is surjective, where $\S_{(\pi/2,\pi/2)}^*\coloneqq\{ g \in \S^* \colon \int g \d\nu=0 \}$ and $\S_{\gamma}^* \coloneqq \S^*$ for $\gamma\not=(\pi/2,\pi/2)$. For $\gamma=(\alpha, \beta) \in [0,{\pi/2}]^{2}\setminus\{(\pi/2,\pi/2)\}$ we have that $\Delta_{\nu,\Lambda}$ is also injective and its inverse $\Delta_{\nu,\Lambda}^{-1} \colon \S_{\gamma}^* \rightarrow \D^{*}_{\gamma}( \Delta_{\nu,\Lambda} )$ has the following kernel representation:
	\begin{align*}
	\Delta_{\nu,\Lambda}^{-1}g \colon x\mapsto \int K_{\alpha,\beta}(x,y) g(y) \d \nu(y)
	\end{align*}
with continuous kernel $K = K_{\alpha,\beta}$ given, for $\alpha,\beta\in [0,\pi/2)$, by
	\begin{align*}
	K(x,y)\coloneqq 
	 		A_{\alpha,\beta}(1+\tan(\beta)-y)( \tan(\alpha)+x)+\1_{[0,x]}(y)(x-y),
	\end{align*}
 where $A_{\alpha,\beta} \coloneqq {-1}/(1+\tan(\alpha)+\tan(\beta))$, and for the remaining cases by
	\begin{align*}
	K(x,y)\coloneqq 
	\begin{cases}
		-(x+\tan(\alpha) )+\1_{[0,x]}(y)(x-y) &\!\!\mbox{if $\beta=\pi/2$ and $\alpha \in [0,\pi/2)$},\\[0.5em]
		(y-1-\tan(\beta))+\1_{[0,x]}(y)(x-y) &\!\!\mbox{if $\alpha= \pi/2$ and $\beta \in [0,\pi/2)$}.
	\end{cases}
	\end{align*}
 For the Neumann case, when $\alpha=\beta=\pi/2$, the operator
$\Delta_{\nu,\Lambda}$ is not injective with kernel $ \Delta_{\nu,\Lambda}^{-1}(\{0\})= \mathbb{R} \1$. Here, $\mathbb{R} \1$ denotes the set of constant functions.
\end{Lem}

\begin{proof}
Let $\alpha$ and $\beta \in [0,\pi/2)$ and let $ g \in \S^{*}$ be fixed. For $x \in [0,1]$, set
	\begin{align*}
	f(x) = \int K_{\alpha,\beta}(x,y)g(y)\d \nu(y) = b\tan(\alpha)+bx+\int_{[0,x]} (x-y)g(y) \d\nu(y)
	\end{align*}
with
	\begin{align*}
	b = \dfrac{-1}{(1+\tan(\alpha)+\tan(\beta))} \left(\tan(\beta)\int g(y) \d\nu(y)+
	\int (1-y)g(y) \d\nu(y)\right).
	\end{align*}
A direct calculation shows that $f(0) = b \tan(\alpha)$, $\nabla_{\Lambda}f(0) = b$, and
	\begin{align*}
	f(1) = b \tan(\alpha) + b + \int (1-y)g(y)\,d\nu(y) \quad \text{and} \quad
	\nabla_{\Lambda}f(1) = b + \int g(y) \d\nu(y).
	\end{align*}
Hence, $f \in \D^{*}_{(\alpha, \beta)}( \Delta_{\nu,\Lambda} )$ with $\Delta_{\nu, \Lambda}f=g$.

Next, let $\alpha \in [0,\pi/2)$, $\beta=\pi/2$ and $ g \in \S^{*}$ be fixed. For $x \in [0,1]$, set
\begin{align*}
 f(x) &= \int K_{\alpha,\beta}(x,y)g(y)\d \nu(y)\\
 &=-\tan(\alpha)\int g(y)\d\nu(y) -x \int g(y)\d\nu(y) +\int_{[0,x]}(x-y)g(y)\d\nu(y).
\end{align*}
A direct calculation shows that
\begin{align*}
f(0)=-\tan(\alpha)\int g(y)\d\nu(y)
\quad \text{and} \quad
\nabla_{\Lambda} f(0)=-\int g(y)\d\nu(y).
\end{align*}
Hence, 
\begin{align*}
f(0)\cos(\alpha) &+ \nabla_{\Lambda} f(0) \sin(\alpha) \\ &=
 - \tan(\alpha)\int g(y)\d\nu(y) \cos(\alpha)+\sin(\alpha)\int g(y)\d\nu(y)=0.
\end{align*}
Further,
\begin{align*}
  \nabla_{\Lambda}f(1)=\int g(y)\d\nu(y)-\int g(y)\d\nu(y)=0.
\end{align*}
Hence, $f \in \D^{*}_{(\alpha, \beta)}( \Delta_{\nu,\Lambda} )$ with $\Delta_{\nu, \Lambda}f=g$.

The final case we need to consider is when $\alpha=\pi/2$ and $\beta\in [0,\pi/2)$. As above, let $g \in \S^{*}$ be fixed, and for $x \in [0,1]$, set
\begin{align*}
 f(x)&\coloneqq  \int K_{\alpha,\beta}(x,y)g(y)\d \nu(y)\\&= \int (y-1-\tan(\beta)) g(y)\d\nu (y)+\int_{[0,x]}(x-y)g(y)\d\nu(y).
\end{align*}
A direct calculation shows that $\nabla_{\Lambda} f(0)=0$, and 
\begin{align*}
f(1)=-\tan(\beta)\int g(y)\d\nu(y)
\quad \text{and} \quad
  \nabla_{\Lambda}f(1)=\int g(y)\d\nu(y).
\end{align*}
With this at hand we obtain 
\begin{align*}
\cos(\beta)f(1)+
  \sin(\beta) \nabla_{\Lambda}f(1)=-\sin(\beta)\int g(y)\d\nu(y)+
  \sin(\beta) \nabla_{\Lambda}f(1)=0,
\end{align*}
and hence, $f \in \D^{*}_{(\alpha, \beta)}( \Delta_{\nu,\Lambda} )$ with $\Delta_{\nu, \Lambda}f=g$.
\end{proof}

A direct consequence of Lemma \ref{lem:Surjec} is that only under Neumann boundary conditions, one has an eigenfunction with corresponding eigenvalue equal to zero.

Before proceeding, we recall the following abstract facts of linear operators: Let us assume $A\colon\dom(A) \subseteq H \rightarrow H$ is a linear, symmetric and surjective operator on a Hilbert space $H$. From this, one may verify that the annihilator of $\dom(A)$ is trivial, that is $\dom(A)^\perp=\{0\}$, and equivalently, $\dom(A)$ is dense in $H$. Further, one can deduce that $A$ is also self-adjoint: The inclusion $\dom(A) \subseteq \dom(A^*)$ holds by symmetry of $A$, where $A^*$ denotes the adjoint of $A$. For the reverse inclusion, note that, for a fixed $f \in \dom(A^*)$, there exists, by surjectivity of $A$, an element $g \in \dom(A)$ such that $A^*f=Ag$. Using symmetry again, for each $h \in \dom(A)$, we have $\langle f ,A h \rangle = \langle A^* f,h\rangle = \langle A g,h \rangle = \langle g , A h\rangle$ and by surjectivity of $A$ we conclude that $f=g \in \dom(A)$. 
 
We can apply these observations to our situation, namely, for $\gamma \in [0,\pi/2]^2$, consider the setting $H=\S_{\gamma}^w$, $A= \Delta_{\nu, \Lambda}$ and $\dom(A)=\D^w_{\gamma}(\Delta_{\nu, \Lambda} )\cap \S_{\gamma}^w$. When $\gamma \in [0, \pi/2]^2 \setminus \{ (\pi/2,\pi/2)\}$, it follows that $ \Delta_{\nu, \Lambda}$ restricted to $\D^w_{\gamma}(\Delta_{\nu, \Lambda} )$ is a densely defined self-adjoint linear operator on $L^{2}(\nu)$. Under Neumann boundary conditions, namely when $\gamma=(\pi/2,\pi/2)$, we have $\Delta_{\nu, \Lambda}^{-1}(\{0\})= \mathbb{R}\1$ and so $L^2(\nu) = \S_{\gamma}^w \oplus \mathbb{R}\1$. Therefore,
\begin{align*}
  \D^w_{\gamma}(\Delta_{\nu, \Lambda} )=\{f +a \colon f\in \D^w_{\gamma}(\Delta_{\nu, \Lambda} )\cap \S_{\gamma}^w \; \text{and} \; a \in \mathbb{R}\}
\end{align*}
is dense in $L^2(\nu) $. Using this observation, it follows that $\Delta_{\nu, \Lambda} $ restricted to $\D^w_{\gamma}(\Delta_{\nu, \Lambda} )$ is a densely defined self-adjoint linear operator on $L^2(\nu)$. The following proposition summarises these observations.

\begin{proposition}\label{prop:adjoint}
For $\gamma \in [0,{\pi/2}]^{2}$, the densely defined operator \mbox{$\Delta_{\nu,\Lambda} \colon L^{2}(\nu)\to L^{2}(\nu)$} with domain $\D^{w}_{\gamma}( \Delta_{\nu,\Lambda})$ is self-adjoint, non-positive and, in particular, closed.
\end{proposition}

\begin{Cor} If $\gamma \in [0,\pi/2]^2\setminus\{ (\pi/2, \pi/2)\}$, then $R_0\coloneqq -\Delta_{\nu,\Lambda}^{-1} \colon \S_{\gamma}^w \to \S_{\gamma}^w$ is compact and self-adjoint.
\end{Cor}

\begin{proof}
Lemma \ref{lem:Surjec} shows that $R_0$ is a Hilbert-Schmidt operator with continuous (bounded) kernel, and is therefore compact. Further, the symmetry of $\Delta_{\nu, \Lambda}$ in tandem with the fact that $R_0$ is bounded, implies that $R_{0}$ is self-adjoint.
\end{proof}

\begin{Cor}\label{cor:Spectral}
Let $\gamma \in [0,\pi/2]^2$ be fixed. The operator $\Delta_{\nu,\Lambda}$ with domain $\D_{\gamma}^{w}(\Delta_{\nu,\Lambda} )$ gives rise to an orthonormal (possibly finite) basis of eigenfunctions with eigenvalues $\lambda_{n}\leq 0$. If $L^{2}(\nu)$ is not finite dimensional, then we have a countable number of eigenvalues with $\lim_{n \rightarrow \infty} -\lambda_n=\infty$, in particular, $\Delta_{\nu,\Lambda}$ is an unbounded operator. On the other hand, if $L^{2}(\nu)$ is finite dimensional, then $\Delta_{\nu,\Lambda}$ has a finite number of eigenfunctions and is bounded.
\end{Cor}

\begin{proof}
For the case $\gamma=(\alpha,\beta) \in [0,\pi/2]^2\setminus\{ (\pi/2, \pi/2)\}$, if $f \in \D^w_{\gamma}\left(\Delta_{\nu, \Lambda} \right)$ is an eigenfunction of $\Delta_{\nu, \Lambda}$ with corresponding eigenvalue $\lambda<0$, then applying Lemma \ref{lem:Surjec} gives, $\Delta_{\nu, \Lambda}f=\lambda f$ if and only if $\lambda^{-1} f = \Delta_{\nu, \Lambda}^{-1}f$. With this at hand, the required result follows from an application of the spectral theorem for self-adjoint compact operators.

For the case $\gamma=(\alpha,\beta)=(\pi/2,\pi/2)$, we consider the resolvent operator \mbox{$R^{\lambda}_{\nu,\Lambda}\coloneqq(\lambda I -\Delta_{\nu,\Lambda})^{-1}$} with $\lambda>0$. From the integral representation of the resolvent operator $R^{\lambda}_{\nu,\Lambda}\coloneqq(\lambda I -\Delta_{\nu,\Lambda})^{-1}$ with domain $\D^w_{\gamma}\left(\Delta_{\nu, \Lambda} \right)$ given in \cite[\textsection\,1.2]{MR0314125}, one may conclude that $R^{\lambda}_{\nu,\Lambda}$ is compact and self-adjoint, see \cite[Theorem 1, p. 251]{MR574035}. This in tandem with the spectral theorem for self-adjoint compact operators yields the required result.
\end{proof}

\begin{Lem}\label{lem:densecm}
For $\gamma=(\alpha,\beta )\in [0,{\pi/2}]^{2} $, set
	\begin{align*}
	C_{\nu,\Lambda}^{\gamma} 
	\coloneqq \begin{cases}
	\left\{ f \in C_{\nu,\Lambda} \colon f(0)=0, \ f(1)=0\right\}& \text{ if } \alpha=\beta=0, \\[0.5em]
	\left\{ f \in C_{\nu,\Lambda} \colon f(0)=0\right\}& \text{if} \; \alpha=0 \; \text{and} \; \beta \in (0,\pi/2],\\[0.5em]
	\left\{ f \in C_{\nu,\Lambda} \colon f(1)=0\right\} &\text{if} \; \alpha \in (0,\pi/2] \; \text{and} \; \beta=0, \\[0.5em]
	 C_{\nu,\Lambda}& \text{if} \; \alpha \; \text{and} \; \beta \in (0,\pi/2].
	\end{cases}
	\end{align*}
The set $\D^s_{\gamma}(\Delta_{\nu, \, \Lambda})$ is dense in $(C_{\nu,\Lambda}^{\gamma}, \lVert \, \cdot \, \rVert_{\infty})$.
\end{Lem}

\begin{proof} 
This result can be found in \cite[Behauptung 2.4]{MR0314125} without a detailed proof. We sketch a proof of this result for the case $\gamma=(\alpha,\beta) \in (0,\pi/2) \times [0,\pi/2)$; the other cases follow analogously, using the kernel representation of $R_0$ for the case $\gamma \in [0,\pi/2]^2\setminus \{(\pi/2,\pi/2)$ as stated in Lemma \ref{lem:Surjec}, and the kernel representation of the resolvent operator $R_{\nu,\Lambda}^{\lambda}$, for some $\lambda>0$, for the case $\gamma=(\pi/2,\pi/2)$, see for instance \cite[Behauptung 2.4]{MR0314125}.
Let $C_{\nu,\Lambda}'$ denote the dual space of $C_{\nu,\Lambda}$ and fix $\Phi \in C_{\nu,\Lambda}' $ such that, for all $f \in \D^s_{\gamma}(\Delta_{\nu, \, \Lambda})$,
	\begin{align*}
	\Phi(f)=\int f(x) \d\phi(x)=0,
	\end{align*}
where $\phi$ denotes the signed distribution function representing $\Phi$. By definition $\phi$ is of bounded variation and local constant on the complement of $\supp(\nu)$. Letting $A_{\alpha, \beta}$ be as in Lemma \ref{lem:Surjec}, and setting $E = -(1+\tan(\beta))A_{\alpha,\beta}$, from the proof of Lemma \ref{lem:Surjec}, Fubini's theorem and integration by parts, for all $g \in C_{\nu,\Lambda}^{\gamma}$, we have that
	\begin{align} 
	\begin{aligned}
	\int (E + s\,A_{\alpha,\beta})g(s) \d\nu(s)& \int (\tan(\alpha)+x) \d\phi(x)\\
	&= \int \int_{[0, x]} (x-s)g(s) \d \nu (s) \d \phi(x) \label{eq:dual1}\\
	&= \int \left(\phi(1)(1-s)- \int_{[s,1]}\phi(s) \d s \right)g(s) \d\nu(s).
	\end{aligned}
	\end{align}
Further, by a second application of integration by parts,
	\begin{align*}
	\int (\tan(\alpha)+x) \d\phi(x)=\tan(\alpha)\left( \phi(1)-\phi(0)\right)+\phi(1)-\int \phi(s) \d s\eqqcolon B_{\phi}.
	\end{align*}
Combining these identities gives
	\begin{align*}
	\int \left(E\,B_{\phi} +s\,A_{\alpha,\beta}\,B_{\phi} -\phi(1)(1-s)+ \int_{[s,1]}\phi(s) \d s\right)g(s) \d\nu(s)=0.
	\end{align*}
If we consider $g(s) = E\,B_{\phi} +s\,A_{\alpha,\beta}\,B_{\phi} -\phi(1)(1-s)+ \int_{[s,1]}\phi(x) \d x \in \mathcal{C}_{\nu,\Lambda}^{\gamma}$ it follows, that for all $s \in [0,1]$
	\begin{align*}
	\int_{[s,1]}\phi(x) \d x=\phi(1)-E\,B_{\phi} -s( A_{\alpha,\beta}\,B_{\phi}+\phi(1)),
	\end{align*}
which is only possible if $\phi(s)=A_{\alpha,\beta} B_{\phi}+\phi(1)$, for all $s \in [0,1]$. Therefore, $\Phi$ is a Dirac measure in $\{0\}$ with weight $\phi(0)$, and by \eqref{eq:dual1}, for all $f \in \mathcal{C}_{\nu,\Lambda}^{\gamma}$,
	\begin{align*}
	\phi(0) \tan(\alpha) \int\left(E+s \,A_{\alpha,\beta} \right)f(s) \d \nu(s)=0.
	\end{align*}
For the particular choice $f \in \mathcal{C}_{\nu,\Lambda}^{\gamma}$ given by $f \colon s\mapsto E+A_{\alpha,\beta}\,s$ the above integral is positive and hence $\phi(0)=0$. Consequently, $\Phi=0$, and since the annihilator of $\D^s_{\gamma}(\Delta_{\nu, \Lambda})$ is trivial, this implies that $\D^s_{\gamma}(\Delta_{\nu, \, \Lambda})$ is dense in $\mathcal{C}_{\nu,\Lambda}^{\gamma}$.
\end{proof}
 
\subsection{Generalised Kre\u{\i}n--Feller operators and transformations of measure spaces}

In this section we return to the general setting where $\mu$ is an arbitrary atomless Borel probability measure with $\supp(\mu) \subseteq [0,1]$ and $\nu$ is embedded in $\mu$. The following two lemmas play a crucial role in the proof of our main result, Theorem \ref{thm:laplace_backward}, which follows directly after.

\begin{Lem}\label{lem:identity}
The function $\check{F}^{-1}_{\mu} \circ F_{\mu}$ equals the identity $\nu$-almost everywhere.
\end{Lem}

\begin{proof}
For $x\in[0,1]$, we have that $\check{F}^{-1}_{\mu}( F_{\mu}(x) )\neq x$ if and only if there exists $\varepsilon>0$ with $F_{\mu}(x-\varepsilon)=F_{\mu}(x)$. This means, if $\check{F}^{-1}_{\mu}( F_{\mu}(x) )\neq x$, then $x$ belongs to an interval of constancy for $F_\mu$. This in tandem with our hypothesis that $\nu$ is embedded in $\mu$ implies that the countable union of the closure of these intervals has $\nu$-measure zero.
\end{proof}

\begin{Lem}\label{lem: isomorphism}
For $*\in \{s,w\}$, the mapping $\varphi \colon \S^{*}(\nu,\mu) \to \S^{*}(\nu\circ {F}^{-1}_{\mu} ,\Lambda)$ defined by \mbox{$\varphi(f) \coloneqq f \circ \check{F}^{-1}_{\mu}$,} is an isometric isomorphism with inverse $\varphi^{-1}(f) = f \circ F_{\mu}$.
\end{Lem}

\begin{proof}
This is a consequence of Lemma \ref{lem:identity} together with the push-forward formula for measures in the weak case and the definition of $C_{\nu,\mu}$ in the strong case.
\end{proof}

The following example demonstrates that Lemmas \ref{lem:identity} and  \ref{lem: isomorphism} do not necessarily hold if the assumption that $\nu$ is embedded in $\mu$ is removed.

\begin{Exa} 
Define the absolutely continuous measure $\mu\coloneqq 2 \Lambda\vert_{[1/2,1]}$ and let $\nu$ denote the Dirac point mass at $1/2$. Note that $\nu$ is not embedded in $\mu$, but that $\supp(\nu) \subseteq \supp(\mu)$ and
$\nu(\{ x\in [0,1] \colon \check{F}^{-1}_{\mu} \circ F_{\mu}(x)\neq x \})=\nu((0,1/2])\neq 0$. 
Further, in the weak setting, since $\nu \circ \check{F}^{-1}_{\mu}$ exhibits an atom at zero, but $f \circ \check{F}_{\mu}^{-1}$ is not defined at zero, the functional $\varphi$ in Lemma \ref{lem: isomorphism} would not be well-defined.
\end{Exa}

\begin{theorem}\label{thm:laplace_backward}
For $\gamma \in [0,{\pi/2}]^{2}$ and $* \in \{ s, w\}$, we have that $\Delta_{\nu \circ F_{\mu}^{-1}, \Lambda }\circ\varphi=\varphi\circ \Delta_{\nu,\mu }$ and $\varphi ( \D^*_{\gamma}(\Delta_{\nu,\mu}))=\D^*_{\gamma}(\Delta_{\nu \circ F_{\mu}^{-1},\Lambda})$.
\end{theorem}

\begin{proof} 
Since $\nu$ is embedded in $\mu$, we have $\nu(F_{\mu}^{-1}(\{0,1\}))=0$. 
If $f \in \D^*_{\gamma}(\Delta_{\nu,\mu})$, then, for all $x \in \supp(\mu)$,
	\begin{align*}
	f(x)=a+bF_{\mu}(x)+
	\int \1_{[0, x]}(y)(F_{\mu}(x)-F_{\mu}(y) ) \Delta_{\nu, \mu } f(y) \d \nu(y),
	\end{align*}
where $a=f(0)$ and $b= \nabla_{\mu} f(0)$. Using Lemma \ref{lem:identity} and replacing $x$ with $\check{F}_{\mu}^{-1}(x)$ gives
	\begin{align*}
	&\hspace*{-1cm}f(\check{F}_{\mu}^{-1}(x))-a-bF_{\mu}(\check{F}_{\mu}^{-1}(x))\\
	&=
	\int \1_{[0, \check{F}_{\mu}^{-1}(x)]}(y)(F_{\mu}(\check{F}_{\mu}^{-1}(x))-F_{\mu}(y) ) \Delta_{\nu, \mu }(f )(y) \d \nu(y) \\
	&=  \int \1_{[0, \check{F}_{\mu}^{-1}(x)]}(y)(x-F_{\mu}(y) ) \Delta_{\nu, \mu }(f )(\check{F}_{\mu}^{-1}(F_{\mu}(y))) \d \nu(y) \\
	&=  \int \1_{[0, x]}(F_{\mu}(y))(x-F_{\mu}(y) ) \Delta_{\nu, \mu }(f )(\check{F}_{\mu}^{-1}(F_{\mu}(y))) \d \nu(y) \\
	&= \int \1_{[0, x]}(y)(x-y ) \Delta_{\nu, \mu }(f ) \circ \check{F}_{\mu}^{-1}(y) \d (\nu \circ F_{\mu}^{-1})(y).
	\end{align*}
This implies that $f(0) = a = f ( \check{F}^{-1}_{\mu}(0) )$, $\nabla_{\mu} f(0) = b=\nabla_{\mu}f ( \check{F}^{-1}_{\mu}(0) )$, \mbox{$f \circ \check{F}_{\mu}^{-1} \in \D^*(\Delta_{\nu \circ F_{\mu}^{-1}, \, \Lambda})$} and
\begin{align*}
\Delta_{\nu \circ {F}_{\mu}^{-1}, \, \Lambda} (f \circ \check{F}_{\mu}^{-1}) = \Delta_{\nu, \, \mu}f\circ \check{F}^{-1}_{\mu}.
\end{align*}
On the other hand, if $g \in \D^*_{\gamma}(\Delta_{\nu \circ F_{\mu}^{-1}, \, \Lambda})$, then by Lemma \ref{lem: isomorphism}, there exists an $f \in \S^{*}(\nu,\mu)$ such that 
$f \circ \check{F}_{\mu}^{-1} = g$, and so, if $x \in \supp(\mu)$ with $F_{\mu}(x-\varepsilon)<F_{\mu}(x)$ for all $\varepsilon>0$, then
	\begin{align}\label{eq:Thm_3_8_star}
	\begin{aligned}
	&\hspace{-2mm}f(x) - c - dF_{\mu}(x) =f(\check{F}_{\mu}^{-1}( F_{\mu}( x)) ) - c - dF_{\mu}(x) \\
	&= \int \1_{[0, F_{\mu}(x)]}(y)(F_{\mu}(x)-y ) \Delta_{\nu \circ F_{\mu}^{-1}, \Lambda}(f \circ \check{F}_{\mu}^{-1})(y) \d (\nu \circ F_{\mu}^{-1})(y)\\
	&= \int \1_{[0, F_{\mu}(x)]}(F_{\mu}(y))(F_{\mu}(x)-F_{\mu}(y) ) \Delta_{\nu \circ F_{\mu}^{-1}, \Lambda}(f \circ \check{F}_{\mu}^{-1})
 \circ F_{\mu}(y) \d \nu (y)\\
	&=  \int \1_{[0, x]}(y)(F_{\mu}(x)-F_{\mu}(y) ) \Delta_{\nu \circ F_{\mu}^{-1}, \Lambda}(f \circ \check{F}_{\mu}^{-1}) \circ F_{\mu}(y) \d \nu (y),
	\end{aligned}
	\end{align}
where $c=g(0) = f (\check{F}^{-1}_{\mu}(0 ))$ and $d = \nabla_{\Lambda}g(0) = \nabla_{\Lambda} f (\check{F}^{-1}_{\mu} (0 ) )$. However, if $x \in \supp(\mu)$ with $F_{\mu}(x-\varepsilon)=F_{\mu}(x)$ for some $\varepsilon>0$, then this implies that $x$ lies in an interval of constancy of $F_{\mu}$. Since $\check{F}_{\mu}^{-1}( [0,1])$ is equal to the set
	\begin{align*}
	(\{0\} \cup \supp(\mu)) \setminus \{ x \in [0,1] \colon x \; \text{  right endpoint of an interval of constancy of} \; F_{\mu}\},
	\end{align*}
we may modify $f$ so that $f$ is constant on each interval of constancy of $F_{\mu}$ while retaining the chain of equalities given in \eqref{eq:Thm_3_8_star}. In other words, given a $g \in \D^*_{\gamma}(\Delta_{\nu \circ F_{\mu}^{-1}, \, \Lambda})$, and letting $f \in \S^{*}(\nu,\mu)$ with $f \circ \check{F}_{\mu}^{-1} = g$, we have $f(\check{F}^{-1}_{\mu}(0 )) = c=f(0)$, $\nabla_{\Lambda} f (\check{F}^{-1}_{\mu} (0 ) ) = d=\nabla_{\mu}f(0)$, $f \in \D^{*}(\Delta_{\nu, \mu})$ and $\Delta_{\nu \circ F_{\mu}^{-1}, \Lambda}(f \circ \check{F}_{\mu}^{-1} ) \circ F_{\mu} = \Delta_{\nu,\mu} ( f )$.

All that remains is to verify the boundary conditions. The above with Lemma \ref{lem:identity} yields
	\begin{align*}
	\nabla_{\mu}f(1) &= \nabla_{\mu}f(0)+\int \Delta_{\nu,\mu}f(y) \d \nu(y)\\
	&= \nabla_{\Lambda}f(0 ) +\int \Delta_{\nu,\mu}f ( \check{F}_{\mu}^{-1}(y) ) \d \nu \circ F_{\mu}^{-1}(y) \\
	&= \nabla_{\Lambda}f(0 ) +\int \Delta_{\nu \circ {F}_{\mu}^{-1}, \, \Lambda} (f \circ \check{F}_{\mu}^{-1}) \d \nu \circ F_{\mu}^{-1}
	= \nabla_{\Lambda} (f \circ \check{F}_{\mu}^{-1} )(1),
	\end{align*}
and similarly that
	\begin{align*}
	&f(1) = f(0) + \nabla_{\mu}f(0)+\int (1 - F_{\mu}(y)) \Delta_{\nu,\mu}f(y) \d \nu(y)\\
	&= f(\check{F}_{\mu}^{-1}(0)) + \nabla_{\Lambda}f(\check{F}_{\mu}^{-1}(0))+\int (1-F_{\mu}(y))\Delta_{\nu,\mu}f ( \check{F}_{\mu}^{-1} \circ F_{\mu} (y) ) \d \nu(y) \\
	&= f(\check{F}_{\mu}^{-1}(0)) + \nabla_{\Lambda}f(\check{F}_{\mu}^{-1}(0))+\int (1-y)\Delta_{\nu,\mu}f ( \check{F}_{\mu}^{-1} (y) ) \d \nu\circ F_{\mu}^{-1}(y) \\
	&= f(\check{F}_{\mu}^{-1}(0)) + \nabla_{\Lambda}f(\check{F}_{\mu}^{-1}(0))+\int (1-y)\Delta_{\nu \circ F_{\mu}^{-1},\Lambda} (f \circ \check{F}_{\mu}^{-1}) (y) ) \d \nu\circ F_{\mu}^{-1}(y)\\&
	= f \circ \check{F}_{\mu}^{-1}(1).\hspace*{8.7cm}\square
	\end{align*}
\end{proof}

\begin{Cor}\label{Cor:cor_3_9}
For $\gamma=(\alpha,\beta )\in [0,{\pi/2}]^{2} $, set
	\begin{align*}
	C_{\nu,\mu}^{\gamma} 
	\coloneqq \begin{cases}
	\left\{ f \in C_{\nu,\mu} \colon f(0)=0, \ f(1)=0\right\}& \text{if} \; \alpha=\beta=0, \\[0.5em]
	\left\{ f \in C_{\nu,\mu} \colon f(0)=0\right\}& \text{if} \; \alpha=0 \; \text{and} \; \beta \in (0,\pi/2],\\[0.5em]
	\left\{ f \in C_{\nu,\mu} \colon f(1)=0\right\} &\text{if} \; \alpha \in (0,\pi/2] \; \text{and} \; \beta=0, \\[0.5em]
	 C_{\nu,\mu} & \text{if} \; \alpha \; \text{and} \; \beta \in (0,\pi/2].
	\end{cases}
	\end{align*}
The set $\D^s_{\gamma}(\Delta_{\nu, \, \mu})$ is dense in $(C_{\nu,\mu}^{\gamma}, \lVert \, \cdot \, \rVert_{\infty})$.
\end{Cor}

\begin{proof}
This follows from Lemma \ref{lem:densecm}, Lemma \ref{lem: isomorphism} and Theorem \ref{thm:laplace_backward}.
\end{proof}

\begin{Cor} 
For each $\gamma \in [0,{\pi/2}]^{2}$, 
the operator $\Delta_{\nu, \, \mu}$ with domain $\D_{\gamma}^w(\Delta_{\nu,\mu})$ is densely defined and self-adjoint.
\end{Cor}

\begin{proof} 
Denseness follows by combining the results Proposition \ref{prop:adjoint}, Lemma \ref{lem: isomorphism} and Theorem \ref{thm:laplace_backward}. To show that $\Delta_{\nu, \, \mu}$ with domain $\D_{\gamma}^w(\Delta_{\nu,\mu})$ is self-adjoint, let us assume the setting of Proposition \ref{prop:adjoint}. For $f \in \D_{\gamma}^w(\Delta_{\nu,\mu})$ by Lemma \ref{lem:identity} and Theorem \ref{thm:laplace_backward}, 
	\begin{align*}
	g\mapsto
	\langle f ,\Delta_{\nu,\mu} g \rangle_{L^2(\nu)}
	=\langle f\circ \check{F}^{-1}_{\mu} , \Delta_{\nu \circ \check{F}^{-1}_{\mu},\Lambda} (g \circ \check{F}^{-1}_{\mu}) \rangle_{L^2(\nu \circ F_{\mu}^{-1})}
	\end{align*}
defines a continuous linear functional on $\D_{\gamma}^w(\Delta_{\nu,\mu})$. Combining Proposition \ref{prop:adjoint} and Theorem \ref{thm:laplace_backward} we may deduce that $ f \circ \check{F}^{-1}_{\mu} \in \dom ( \Delta_{\nu \circ \check{F}^{-1}_{\mu},\Lambda}^{*} ) =\D_{\gamma}^w( \Delta_{\nu \circ \check{F}^{-1}_{\mu},\Lambda})$, and consequently that $f \in \D_{\gamma}^w(\Delta_{\nu,\mu})$.
\end{proof}

By taking into account that $\Delta_{\nu,\mu}$ with domain $\D_{\gamma}^{*}(\Delta_{\nu,\mu})$ has a probabilistic interpretation as infinitesimal generator of a Markov process, one may also obtain the above corollary via an application of \cite[Theorem 10.13]{DynkinI_II}. 

Corollary \ref{cor:Spectral} implies that $\Delta_{\nu \circ {F}_{\mu}^{-1}, \, \Lambda}$ with domain $\D_{\gamma}^w \left( \Delta_{\nu, \Lambda} \right)$ gives rise to an orthonormal basis of eigenfunctions with non-positive eigenvalues. In the case that $\supp(\nu)$ is infinite, we have $(\lambda_n)_{n \in \mathbb{N}}$ with $\lim_{n \to \infty}-\lambda_n = \infty$; otherwise there are only finitely many eigenvalues. Using the one-to-one correspondence established in Theorem \ref{thm:laplace_backward} to relate the spectral properties of $\Delta_{\nu,\mu}$ with those of $\Delta_{\nu \circ {F}_{\mu}^{-1}, \, \Lambda}$, we obtain the following.

\begin{Cor}\label{cor:Spec}
For fixed $ \gamma \in [0,{\pi/2}]^{2}$, the operators $\Delta_{\nu \circ {F}_{\mu}^{-1}, \, \Lambda} $ with domain $\D^w_{\gamma}( \Delta_{\nu \circ F_{\mu}^{-1},\Lambda} )$ and $ \Delta_{\nu,\mu} $ with domain $\D^w_{\gamma}( \Delta_{\nu,\mu})$ have the same eigenvalues $(\lambda_n)_{n \in \mathbb{N}}$. Further, if $f$ is an eigenfunction of $ \Delta_{\nu,\mu} $, then $f \circ \check{F}^{-1}_{\mu}$ is an eigenfunction of $ \Delta_{\nu \circ {F}_{\mu}^{-1}, \, \Lambda} $, and if $f$ is an eigenfunction of $\Delta_{\nu \circ {F}_{\mu}^{-1}, \, \Lambda}$, then $f \circ F_{\mu}$ is an eigenfunction of $ \Delta_{\nu,\mu}$. In particular, if $(f_n)_{n \in \mathbb{N}}$ denotes the orthonormal basis consisting of eigenfunctions of $ \Delta_{\nu \circ {F}_{\mu}^{-1}, \, \Lambda}$, then $(f_n \circ F_{\mu})_{n \in \mathbb{N}}$ forms an orthonormal basis comprising of eigenfunctions of $\Delta_{\nu,\mu}$.
\end{Cor}
 
Corollary \ref{cor:Spec} can be seen as a generalisation of the results of \cite{KSW16} where $\Delta_{\mu,\mu}$ was considered.

\begin{Cor}\label{cor:Reso} For $\gamma \in[0,{\pi/2}]^{2}$ and $\lambda>0$, letting $R_{\nu , \mu}^{\lambda} = (\lambda - \Delta_{\nu, \mu})^{-1}$ denote the resolvent operator of $ \Delta_{\nu, \mu}$ with domain ${ \D^*_{\gamma}(\Delta_{\nu,\mu})}$, for all $f \in \S^*$,
	\begin{align}\label{eq:resolvent_equation}
	R_{\nu,\mu}^{\lambda}( f ) \circ \check{F}_{\mu}^{-1}=R_{\nu \circ F_{\mu}^{-1},\Lambda}^{\lambda}(f \circ \check{F}_{\mu}^{-1}).
	\end{align}
In particular, we have that
	\begin{align*}
	\Vert R_{\nu \circ F_{\mu}^{-1},\Lambda}^{\lambda} \Vert_{L^2(\nu \circ F_{\mu}^{-1} ) }= \Vert R_{\nu,\mu}^{\lambda} \Vert_{L^2 (\nu)}
	\quad \text{and} \quad
	\Vert R^{\lambda}_{\nu \circ F_{\mu}^{-1},\Lambda} \Vert_{C_{\nu \circ F_{\mu}^{-1},\Lambda}}
	= \Vert R^{\lambda}_{\nu,\mu}\Vert_{C_{\nu,\mu}}.
	\end{align*} 
\end{Cor}

\begin{proof}
The resolvent $R^{\lambda}_{\nu,\mu}$ is well defined for all $\lambda>0$, as  by Proposition \ref{prop:adjoint} the operator $\Delta_{\nu,\mu}$ is non-positive, and for $f \in \S^*$, we have
	\begin{align*}
	\Delta_{\nu,\mu} (R_{\nu ,\mu}^{\lambda}(f))
	=\lambda R_{\nu, \mu}^{\lambda} (f) - f
	\quad \text{and} \quad
	R_{\nu ,\mu}^{\lambda} (f) \in \D^*_{ \gamma}( \Delta_{\nu,\mu} ).
	\end{align*}
Theorem \ref{thm:laplace_backward} gives 
	\begin{align*}
	\Delta_{\nu\circ F_{\mu}^{-1},\Lambda} (R_{\nu ,\mu}^{\lambda}(f) \circ \check{F}^{-1}_{\mu}) 
	= \Delta_{\nu,\mu} (R_{\nu ,\mu}^{\lambda}(f)) \circ \check{F}^{-1}_{\mu} 
	=\lambda R_{\nu,\mu}^{\lambda}(f)\circ \check{F}^{-1}_{\mu} -f \circ \check{F}^{-1}_{\mu},
	\end{align*}
proving the first part. This in tandem with Lemma \ref{lem: isomorphism} yields \eqref{eq:resolvent_equation}.
\end{proof}

\begin{Cor}\label{cor:Feller}
For $\gamma \in[0,{\pi/2}]^{2}$, the operator $\Delta_{\mu,\nu}$ with domain $ \D_{\gamma}^s(\Delta_{\nu,\mu})$ is an infinitesimal generator of a strongly continuous semi-group of contractions on $C_{\nu,\mu}^{\gamma}$.
\end{Cor}

\begin{proof}
It is well-known that the result hold for the classical Kre\u{\i}n--Feller operator, see \mbox{\cite[Behauptung 4.1]{MR0314125}}. This in tandem with the Yosida-Hille Theorem, see for example \cite[Theorem 1.12]{ ma1992introduction}, and Corollary \ref{cor:Reso}, yields 
	\begin{align*} 
	\Vert R^{\lambda}_{\nu \circ F_{\mu}^{-1},\Lambda} \Vert_{C^{\gamma}_{\nu \circ F_{\mu}^{-1},\Lambda}}
	= \Vert R^{\lambda}_{\nu,\mu} \Vert_{C_{\nu,\mu}^{\gamma}} \leq 1/\lambda ,
	\end{align*}
for all $\lambda>0$. By Corollary \ref{Cor:cor_3_9}, the operator $\Delta_{\mu,\nu}$ with domain $ \D_{\gamma}^s(\Delta_{\nu,\mu})$ is densely defined, and since convergence in $\lVert \, \cdot \, \rVert$ implies convergence in $\lVert \, \cdot \, \rVert_{L^{2}(\nu)}$, Proposition \ref{prop:adjoint} implies that $\Delta_{\mu,\nu}$ with domain $ \D_{\gamma}^s(\Delta_{\nu,\mu})$ is closed. Thus, a second application of the Yosida-Hille Theorem yields the required result . 
\end{proof}

\section{Applications}

In this last section, some applications to the spectral asymptotics of the Kre\u{\i}n-Feller operators and to the associated gap diffusion are given. 
\subsection{Spectral asymptotics}\label{sec:SpectralAsymp}

In this section, let us begin to link Theorem \ref{thm:laplace_backward} to the result of \cite{KN22A}. For this we require the following.  For $\nu$ and $\mu$ satisfying our standing assumptions, we define the \textsl{upper spectral dimension} $\overline{s}_{\nu,\mu}$ of the generalised Kre\u{\i}n--Feller operator $\Delta_{\nu,\mu}$ by
    \begin{align*}
        \overline{s}_{\nu,\mu}\coloneqq \limsup_{x\to \infty} \frac{\ln(N^\gamma_{\nu,\mu}(x))}{x},
    \end{align*}
    where $N^\gamma_{\nu,\mu}(x)$ denotes the number of eigenvalues of $-\Delta_{\mu,\nu}$  with domain $ \D_{\gamma}^w(\Delta_{\nu,\mu})$    not exceeding $x\in \R$. 
For $q>0$ and $\eta$ a compactly supported measure on $\mathbb{R}$, we denote by 
    \begin{align*}
        \beta_\eta(q) \coloneqq \limsup_{n\to\infty}\frac{1}{n\ln2}\ln \left( \sum_{k\in \mathbb{Z}}\left(\eta\left((k2^{-n},(k+1)2^{-n}]\right)\right)^q\right)
    \end{align*}
the \textsl{$L^q$-spectrum} of $\eta$ and let $q_\eta \coloneqq \inf \{ q>0 : \beta_\eta(q) < q \}$.  Further, we let $\overline{\dim}_M(A)$ denote the \textsl{upper Minkowski dimension} of a bounded set $A\subseteq\mathbb{R}$.

\begin{theorem}\label{thm:USD}
    For $\nu$ and $\mu$ satisfying our standing assumptions, we have that
        \begin{align*}
            \overline{s}_{\nu,\mu} = q_{\nu\circ F_\mu^{-1}}\leq  \frac{\overline{\dim}_M(\supp(\nu\circ F_\mu ^{-1}))}{{1+\overline{\dim}_M(\supp(\nu\circ F_\mu ^{-1}))}}.
        \end{align*}
\end{theorem}

\begin{proof}
    This follows from Theorem \ref{thm:laplace_backward} in combination with \cite[Theorem~1.1 and Corollary.~1.7]{KN22A}.
\end{proof}

Next, we review the asymptotic growth rate of the eigenvalue counting function of 
$-\Delta_{\nu,\mu}$ for certain classes of self-similar measures. With the help of Theorem \ref{thm:laplace_backward}, we show how the results in \cite{MR2787628} can be deduced from \cite{MR1328700}, and hence, how the asymptotic 
results of \cite{Fr05} would follow from \cite{Fu87}. Let us begin by giving our standing assumptions on the self-similar measures. These assumptions are less restrictive than (A.1)~--~(A.4) of \cite{Fr05}, in the sense that our contractions are allowed to be both order-preserving and order-reversing.

\begin{assumptions}\label{ass:StrongAssumption}
Let $M \geq 2$ be an integer and let 
\begin{align*}
\Phi = (S_i \colon [0,1]\to [0,1] \colon i \in \{1,\dots,M\})
\end{align*}
denote a family of affine contractions fulfilling the open set condition (OSC) with feasible open set $(0, 1)$, that is, $S_i((0,1)) \cap S_j((0,1))= \emptyset$, for all $i$ and $j \in \{ 1, \dots M\}$ with $j \neq i$. Let $\nu$ denote the self-similar measure associated with $\Phi$ and with probability weight vector $(p_{1},\dots ,p_{M}) \in (0,1)^M$, namely, the Borel probability measure uniquely determined by
	\begin{align*}
	\nu(A)=\sum_{i=1}^{M} p_i \nu(S_i^{-1}(A))
	\end{align*}
for all $A \in \mathfrak{B}([0,1])$ and where $\mathfrak{B}([0,1])$ denotes the Borel $\sigma$-algebra on $[0, 1]$. For a fixed $(\sigma_{1},\dots ,\sigma_{M})\in (0,1)^{M}$ with $\sum_{i=1}^M \sigma_i \leq 1$, let $\mu$ denote an atomless Borel probability measure with $\supp(\nu) \subseteq \supp(\mu)$ and such that $\mu(S_i(A))=\sigma_i \mu(A)$, for all $A \in \mathfrak{B}([0,1])$ and $i \in \{1,\dots ,M \}$. For $k \in \{ 1, \dots, M \}$, we set $c_{k} \coloneqq - \ln(\sigma_{k} p_{k})$, and distinguish the following two cases.
\begin{enumerate}[leftmargin=*,itemsep=0.5em]
\item \textsl{The arithmetic case}: When the set $\{ c_1, \dots, c_{M} \}$ is a subset of a discrete subgroup of $(\mathbb{R},+)$, that is, there exists a real number $L>0$ such that 
$c_k\in L\,\mathbb{N}$ for all $k \in \{1,\dots,M \}$.
\item \textsl{The non-arithmetic case}: When ${c_k}/{c_l}$ is irrational for some $k$ and $l \in \{ 1, 2, \dots, M\}$.
\end{enumerate}
\end{assumptions}

\begin{theorem}\label{thm:Freiberg}
Assume the setting of Assumption \ref{ass:StrongAssumption}. Let $\gamma \in [0, \pi/2]^{2}$ be fixed and let $u\in (0,1)$ be the unique real number satisfying $\sum_{i=1}^M ( \sigma_i p_i )^{ u}=1$.  Let $( \lambda_n )_{n \in \mathbb{N}}$ denote the increasing sequence of eigenvalues of $-\Delta_{\nu, \mu}$ with domain $\D^{w}_{\gamma}(\Delta_{\nu ,\mu})$ and denote by $N_{\nu,\mu}^{\gamma}$ the associated eigenvalue counting function. In the non-arithmetic case, there exists a positive constant $k$ such that
	\begin{align*}
	\lim_{n \to \infty} n^{-1/u} \lambda_n = k
	\quad \text{and} \quad
	\lim_{x \to \infty} x^{-u} N_{\nu,\mu}^{\gamma}(x) = k^{-u}.
	\end{align*}
In the arithmetic case, there exists an $L$-periodic function $\psi$ bounded and separated from zero, such that
	\begin{align*} 
	\lim_{x \to \infty}
	x^{-u} N_{\nu,\mu}^{\gamma}(x)/ \psi( \ln(x)) = 1.
	\end{align*}
\end{theorem}

\begin{proof}
For $i \in \{1, \dots, M\}$, let $\mathrm{sgn}(S_{i})$ denote the sign of the affine transformation $S_{i}$, that is, if $S_{i}$ is order preserving, then $\mathrm{sgn}(S_{i}) = 1$, and otherwise $\mathrm{sgn}(S_{i}) = -1$. With this at hand, we observe that
	\begin{align*}
	\sigma_i F_{\mu}(x)=\sigma_i \mu([0,x])= \mu(S_i( [0,x] ) )= \mathrm{sgn}{(S_{i})} (F_{\mu}(S_i(x))- F_{\mu}(S_i(0))),
	\end{align*}
for all $x \in [0,1]$ and $i \in \{1,\dots,M\}$. Setting $\tilde{S}_i(x) = \mathrm{sgn}(S_{i}) \sigma_i x+F_{\mu}(S_i(0))$, for $x \in [0,1]$ and $i \in \{1,\dots,M\}$, we obtain $F_{\mu} \circ S_i= \tilde{S}_i \circ F_{\mu}$, and hence, for all $A \in \mathfrak{B}([0,1])$, that
	\begin{align*}
	\sum_{i=1}^M p_i\nu( F_{\mu}^{-1} \circ \tilde{S}_i^{-1}(A)) 
	&=\sum_{i=1}^M p_i\nu( S_i^{-1} \circ F_{\mu}^{-1}(A))
	=\nu ( F_{\mu}^{-1}(A) ). 
	\end{align*}
This shows that $ \nu \circ F_{\mu}^{-1}$ is the unique self-similar measure for the family of contractions $( \tilde{S}_i \colon i \in \{1,\dots, N \} )$ and probability weight vector $(p_{1},\dots ,p_{M})$. Additionally, since $\supp(\nu) \subseteq \supp(\mu)$ and since
	\begin{align*}
	\tilde{S}_i( [0,1] )
	&= [ F_{\mu}(S_i(0)) + (\mathrm{sgn}(S_{i})-1)\sigma_{i}/2, F_{\mu}(S_i(0)) + (\mathrm{sgn}(S_{i})+1)\sigma_{i}/2 ]\\
	&=[\min \{ F_{\mu}(S_i(0)),F_{\mu}(S_i(1)) \}, \max \{ F_{\mu}(S_i(0)),F_{\mu}(S_i(1)) \}],
	\end{align*}
$( \tilde{S}_i \colon i \in \{1,\dots ,M\} )$ satisfies the OSC given that $( {S}_i \colon i \in \{1,\dots ,M\} )$ satisfies the OSC with feasible open set $(0, 1)$. Therefore, we can apply the classical result from Solomyak and Verbitsky \cite{MR1328700} for the spectral asymptotics of $\Delta_{\nu \circ {F}_{\mu}^{-1}, \, \Lambda} $ with domain ${\D^{w}_{(0,0)} (\Delta_{\nu \circ {F}_{\mu}^{-1}, \, \Lambda} )}$. The spectral asymptotics for the Robin case $\gamma \in [0,\pi/2]^2$ follows by well-known estimates of the eigenvalue counting function, see for example \cite[Lemma 3.1]{Fu87}. Combining this with Corollary \ref{cor:Spec} completes the proof.
\end{proof}

\begin{Exa}\label{exa:selfSim}
Fix $N$ and $M \in \mathbb{N}$ with $2 \leq M \leq N$. For $i \in \{1,\dots,N\}$, we let \mbox{$S_{i} \colon [0,1] \to [0, 1]$} be defined by $S_{i}(x) = s_i x+b_i$, where $s_i$ and $b_i$ are non-negative real numbers with $s_i+b_i \leq 1$, and such that $( S_{i} \colon i \in \{1, \dots, N\})$ fulfils the OSC with feasible open set $(0, 1)$. Let $\delta_{N}$ and $\delta_{M}$ be the unique positive real numbers satisfying
	\begin{align*}
	\sum_{i=1}^{N} s_{i}^{\delta_{N}} = 1 \quad \text{and} \quad \sum_{i=1}^{M} s_{i}^{\delta_{M}} = 1.
	\end{align*}
Denote by $\nu$, the unique self-similar measure determined by the family of contractions $(S_{i} \colon i \in \{ 1, \dots ,M\} )$ and probability weight vector 
\begin{align*}
(p_{1},\dots ,p_{M})\coloneqq (s_{1}^{\delta_{M}}, \dots, s_{M}^{\delta_{M}}),
\end{align*}
and let $\mu$ be the unique self-similar measure determined by the family of contractions $(S_{i}\colon i \in \{ 1,\dots ,N \} )$ and probability weight vector $( s_{1}^{\delta_{N}},\dots, s_{N}^{\delta_{N}})$. In which case, $\supp(\nu) \subseteq \supp(\mu)$,
	\begin{align*}
	\nu \circ F_{\mu}^{-1}(A)= \sum_{i=1}^M p_{i}
	\nu(S_i^{-1}\circ F_{\mu}^{-1}(A)),
	\end{align*}
and hence, via an application of Theorem \ref{thm:Freiberg}, we find that the exponential growth rate of the eigenvalue counting function of $-\Delta_{\nu, \mu}$ is 
$\delta_{M}/(\delta_{N}+\delta_{M})$.  Additionally, we note that, according to the Moran--Hutchinson formula
$\overline{\dim}_M(\supp(\mu)) = \delta_{N}$ and $\overline{\dim}_M(\supp(\nu)) = \delta_{M}$.
\end{Exa}

The estimates on the upper spectral dimension of the classical Kre\u{\i}n--Feller operator, as described in \cite[Corollary~1.8]{KN22A}, sheds new light on Example \ref{exa:selfSim}.
Indeed, if $\mu$ is $\delta_\mu$-Ahlfors regular (as in Example \ref{exa:selfSim}), which in the one-dimensional context means that $F_\mu$ is $\delta_\mu$-Hölder continuous, then  
    \begin{align*}
        \overline{\dim}_M(\supp(\nu\circ F_\mu ^{-1}))
        &= \overline{\dim}_M(F_\mu(\supp(\nu))
        \\&\leq \frac{\overline{\dim}_M(\supp(\nu))}{\delta_\mu}=\frac{\overline{\dim}_M(\supp(\nu))}{\overline{\dim}_M(\supp(\mu))}.
    \end{align*}
Here, we have also used that $\delta_\mu = \overline{\dim}_M(\supp(\mu))$ which follows from the defining property of Ahlfors regularity and the equivalent definitions of the upper, and the lower, Minkowski dimension given by centered $\delta$-packings, and $\delta$-covers, respectively.
So, by Theorem \ref{thm:USD}, if $\mu$ is Ahlfors regular, then
    \begin{align*} 
        \overline{s}_{\nu,\mu}\leq \frac{\overline{\dim}_M(\supp(\nu))}{\overline{\dim}_M(\supp(\nu))+\overline{\dim}_M(\supp(\mu))}.
    \end{align*} 
In this respect, fixing  $\mu$ and  $\supp(\nu)$, the upper spectral dimension $\overline{s}_{\nu,\mu}$ is maximal in our example. We emphasise that the upper Minkowski dimension, and not the Hausdorff dimension, is the relevant quantity in this estimate. For further interesting connections of the the $L^q$-spectrum of a compactly supported measures $\nu$ with its quantization, i.\,e.\ its approximations by finitely supported measures, we refer the interested reader to \cite{KNZ23}.

As a third application of our main result, we provide a simple proof of the folklore result that the absolutely continuous part of $\nu$ with respect to $\mu$ is dominating the spectral asymptotic, see \cite{Volkmer05} for an alternative approach in a slightly more general setting. 

\begin{theorem}\label{thm:RayMacKean}
Let $\mu$ and $\eta$ denote two finite atomless  Borel   measures on $\mathbb{R}$ such that we have \mbox{$\supp(\eta) \subseteq \supp(\mu) \subseteq [0, 1]$} and $\eta$ is singular to $\mu$, and let $\sigma\in L^{2}(\mu) $ be  non-negative. Assuming both  $\nu \coloneqq  \sigma^{2} \mu + \eta$ and $\mu$ are probability measures, 
then for $\gamma \in \{ (0,0), (\pi/2, \pi/2)\}$,
	\begin{align*}
	\lim_{x\to\infty}\frac{N_{\nu,\mu}^{\gamma}(x)}{x^{1/2}}=   \pi^{-1} \int  \sigma \d\mu.
	\end{align*}
\end{theorem}

\begin{proof}
For $x \in \mathbb{R}$ and $\gamma \in \{ (0,0), (\pi/2, \pi/2)\}$, by Corollary \ref{cor:Spec},
	\begin{align*}
	N_{\sigma^{2}\mu+\eta,\mu}^{\gamma}(x)=N_{\sigma^{2}\circ\check{F_{\mu}}^{-1}\Lambda+\eta\circ F_{\mu}^{-1},\Lambda}^{\gamma}(x).
	\end{align*}
By our hypothesis there exists a Borel set $A \subseteq [0,1]$ with $\eta(A)=0$ and $\mu( A)=1$. By an application of Lemma \ref{lem:identity},
	\begin{align*}
	\eta \circ F_{\mu}^{-1}((\check{F}_{\mu}^{-1})^{-1} A) = \eta (F_{\mu}^{-1}((\check{F}_{\mu}^{-1})^{-1}A))=\eta((\check{F}_{\mu}^{-1}\circ F_{\mu})^{-1}A)=\eta(A)=0
	\end{align*}
and 
	\begin{align*}
	\Lambda((\check{F}_{\mu}^{-1})^{-1}A)=\mu\circ F_{\mu}^{-1}((\check{F}_{\mu}^{-1})^{-1}A)=\mu((\check{F}_{\mu}^{-1}\circ F_{\mu})^{-1}A)=\mu(A)=1.
	\end{align*}
Hence, the measure $\eta\circ F_{\mu}^{-1}$ is singular to $\Lambda$, and thus, by \cite[Theorem~5.1]{MR0278126}, or alternatively the results of \cite{Volkmer05}, we obtain, for $\gamma \in \{ (0,0), (\pi/2, \pi/2)\}$,
	\[
		\hspace*{5mm}\lim_{x\to\infty}\frac{N_{\sigma^{2} \circ \check{F}_{\mu}^{-1}\Lambda+\eta\circ F_{\mu}^{-1},\Lambda}^{\gamma}(x)}{x^{1/2}}
		= \pi^{-1} \int \left(\sigma^{2} \circ \check{F_{\mu}}^{-1}\right)^{1/2} \d \Lambda
		= \pi^{-1}\int \sigma \d\mu.\hspace*{7mm}\square 
	\]
\end{proof}

\subsection{Gap diffusion}\label{sec:LBM}

In this section we define generalised $\nu$-$\mu$-gap diffusion via a time change of a fixed Brownian motion and a state space transformation given by $\check{F}^{-1}_\mu$.
Making use of known basic properties for gap diffusion and our general approach, we show that the infinitesimal generator of the generalised gap diffusion coincides with a generalised Kre\u{\i}n--Feller operator. Our standing assumptions for this section are as follows.

For a measurable space  $(\Omega, \mathcal{F})$ and for some  $x \in \Omega$, we consider the tuple $\left(\Omega, \mathcal{F}, (\mathcal{F}_t)_{t \geq 0}, ( \theta_t)_{t \geq 0}, \mathbb{P}_x\right)$ where $\mathbb{P}_x$ denotes the probability measure such that, $(\mathcal{F}_t)_{t \geq 0}$ is the right-continuous completed natural filtration and $(\theta_t )_{t \geq 0}$ is the shift-operator. The expectation with respect to $\mathbb{P}_x$ is denoted by $\mathbb{E}_x$, and we call a stochastic process $(B_t )_{t \geq 0}$ a \textsl{Brownian motion} if:
\begin{enumerate}[leftmargin=*,itemsep=0.5em]
\item $\mathbb{P}_x(B_0=x)=1$
\item For $0 \leq s_0< \dots <s_n$ with $n \in \mathbb{N}$, the increments $B_{s_1}-B_{s_0}, \dots ,B_{s_n}-B_{s_{n-1}}$ are stochastically independent.
\item For $s$ and $t \in \mathbb{R}$ with $s > t \geq 0$ the distribution of the increment $B_{s}-B_{t}$ follows a Gaussian distribution with variance $s-t$ and mean $0$.
\end{enumerate} 
Further, we let $(L_{x}^t)_{t \geq 0, x \in \mathbb{R}}$ be the jointly continuous version of the local time of the Brownian motion $(B_t )_{t \geq 0}$, see \cite[Chapter VI]{Revuz2013}, and let $m$ be a probability Borel measure on $[0,1]$.

We define, for $ t \geq 0$, the \textsl{time-change function}
	\begin{align*}
	\Phi_t \coloneqq \int L_x^t \dm(x)
	\end{align*}
and the \textsl{right-continuous inverse} of $\Phi_t$,
	\begin{align*}
	\hat{\Phi}^{-1}_t \coloneqq \inf\left\{s \geq 0 \colon \Phi_s>t\right \}.
	\end{align*}
For $x \in \supp(m)$, the process $(X_t)_{t\geq 0}\coloneqq(B_{\hat{\Phi}^{-1}_t})_{t \geq 0}$, defined with respect to the filtration $(\mathcal{F}_{\hat{\Phi}^{-1}_t})_{t \geq 0}$ and probability measure $\mathbb{P}_x$, is called a {\em gap diffusion} with {\em speed measure} $m$ and {\em starting point} $x$.
 
Now we recall some important facts for gap diffusion relevant for our considerations.

\begin{Fact}[{\cite[Lemma 3.1]{MR3034785}}]For all $t \geq 0$, we have $\mathbb{P}_x$-almost everywhere, that $X_t \in \supp(m)$.\label{fact:SkipfreeP}
\end{Fact}

\begin{Fact}[{\cite[Satz 2.4.1, Satz 2.4.6]{Burkhardt1983}}]
The tuple 
\begin{align*}
((X_t )_{t \geq 0}, (\mathcal{F}_{\hat{\Phi}^{-1}_t})_{t \geq 0}, (\theta_{\hat{\Phi}^{-1}_t} )_{t \geq 0}, ( \mathbb{P}_x )_{x \in \supp(m)})
\end{align*}
defines a {\em Feller process} with state space $E\coloneqq \supp(m)$ equipped with the Euclidean topology $\tau$. Namely, $(X_t )_{t \geq 0}$ is a strong Markov process such that, for all $f$ from the set $C_b (E )$ of bounded continuous function with domain $E$, the function $x \mapsto \mathbb{E}_x[f(X_t) ]$ belongs to $C_b(E )$ and for every $x\in E$ we have \begin{align*} \lim_{ t \searrow 0} \mathbb{E}_x[f(X_t)] =f(x). \end{align*}
\end{Fact}

We now want to make explicit a connection between generalized Kre\u{\i}n--Feller operators and the infinitesimal generator of transformed gap diffusions. To this end, let $\mu$ and $\nu$ denote two Borel probability measures satisfying our standing assumptions given at the start of Section \ref{sec:standing_assumptions}, and let $(X_t)_{ t \geq 0}$ denote a gap diffusion with speed measure $\nu \circ F_{\mu}^{-1}$. We call the stochastic process $(\check{F}_{\mu}^{-1}(X_t))_{t \geq 0}$ a (\textsl{generalised}) \textsl{$\nu$-$\mu$-gap diffusion} with  speed measure $\nu$ and scale measure $\mu$. We will see below that by \cite[Vol.\,I, \textsection\,6, Theorem 10.13]{DynkinI_II} this latter stochastic process is again strong Markov. It should be noted that the fact that ${F}_{\mu}$ is continuous is essential to ensure that $(Y_t)_t$ is a strong Markov process (see also \cite{Ogura_1989} for a counterexample if $F_{\mu}$ is not continuous). 

\begin{definition}
Let $(Z_t)_{t \geq 0}$ denote a Markov process with  state space $(E,\tau)$.  
A function $f \in C_b(E)$ is said to belongs to the domain $\D(A)$ of the \textsl{infinitesimal generator} $A$ of $(Z_t)_{t \geq 0}$ if
	\begin{align*}
	Af(x) = \lim_{t \searrow 0} \; (\mathbb{E}_x[f(Z_t) ]-f(x))/t
	\end{align*}
exists with respect to $\lVert \, \cdot \, \rVert_{\infty}$.
\end{definition}

We are now in a position to compute the infinitesimal generator of a \textsl{$\nu$-$\mu$-gap diffusion}. To this end, let $(X_t)_{ t \geq 0}$ denote a gap diffusion with speed measure $\nu \circ F_{\mu}^{-1}$ and, for $t\geq 0$, set \[Y_t \coloneqq \check{F}^{-1}_{\mu}(X_t).\]
Note, the state space of $X_t$ is given by $E\coloneqq \supp(\nu\circ F_\mu^{-1})=F_{\mu}(\supp(\nu))$ equipped with the subspace topology $\tau$ induced by the Euclidean topology, and the state space of $Y_t$ is equal to $\tilde{E}\coloneqq \check{F}_\mu^{-1}(E)$ equipped with the final topology $\tilde{\tau}$ induced by $\check{F}_\mu^{-1}$, i.\,e.\ $\tilde{\tau}\coloneqq\{O\subseteq \tilde{E}:(\check{F}_\mu^{-1})^{-1}(O)\in \tau\}$.

\begin{Lem}\label{lem:Equi_cont_function}
The map $\check{F}_\mu^{-1} \colon E\to \tilde{E}$ defines a homeomorphism, where $E$ and $\tilde{E}$ are respectively equipped with the topologies $\tau$ and $\tilde{\tau}$.
\end{Lem}

\begin{proof}
Note that $\check{F}_\mu^{-1}$ is strictly increasing and hence bijective, and that by definition of the final topology, also continuous. To see that $\check{F}_\mu^{-1}:E\to \tilde{E}$ is a homeomorphism we apply the universal property of the final topology and the fact that $F_{\mu} \circ \check{F}_{\mu}^{-1}$ is the identity. 
\end{proof}

An application of this lemma yields
  \begin{align*}
  C_{\nu\circ F_{\mu}^{-1},\Lambda}\to C(E): f\mapsto f|_E \mbox{ and } C_{\nu,\mu}\to C(\tilde{E}): f\mapsto f|_{\tilde{E}}
  \end{align*}
are bijective, and thus allows us to identify $C(E)$ with $C_{\nu\circ F_{\mu}^{-1},\Lambda}$ and $C(\tilde{E})$ with $C_{\nu,\mu}$. As an aside, we remark that the topology $\tilde{\tau}$ can also be described as the initial topology on $\tilde{E}$ with respect to $\left \{f|_{\tilde{E}}:f\in C_{\nu,\mu}\right\}$. Further, for all $f\in C(\tilde{E})$, we have $f=f \circ \check{F}_{\mu}^{-1} \circ {F}_{\mu}$.

Letting $A_X$ be the infinitesimal generator of $(X_t)_{t \geq 0}$ on $\D(A_{X})\subseteq C_{\nu \circ F_{\mu}^{-1},\Lambda}$, for $f \in \D(A_Y)\subseteq C_{\nu,\mu}$ and $x \in E$, we have $f \circ \check{F}_{\mu}^{-1} \in C_{\nu \circ F_{\mu}^{-1},\Lambda} $ and therefore it follows
	\begin{align*}
	\lim_{t \searrow 0} \, ( \mathbb{E}_x[ f( \check{F}_{\mu}^{-1}(X_t))-f(\check{F}^{-1}_{\mu}(x)) ] )/t
	&=\lim_{t \searrow 0} \, ( \mathbb{E}_{\check{F}_{\mu}^{-1}(x)}[ f (Y_t) -f(\check{F}^{-1}_{\mu}(x))] )/t\\
	&=\lim_{t \searrow 0} \, ( \mathbb{E}_{\check{F}_{\mu}^{-1}(x)}[ f (Y_t) -f(Y_0)] )/t.
	\end{align*}
As in \cite[Vol.\,I, \textsection\,6, Theorem 10.13 \& 1st Remark]{DynkinI_II}, this implies  that $(Y_t)_t$ is a strong Markov process and
\begin{align*}
f\circ \check{F}_\mu ^{-1} \in \D(A_{X}) \iff f\in \D(A_Y) \; \text{and}	\; A_X(f \circ \check{F}_{\mu}^{-1})=A_Y ( f) \circ \check{F}_{\mu}^{-1}, f\in \D(A_Y).
\end{align*}

From \cite[p.\,49]{Burkhardt1983}, we have $ \D(A_X)\subseteq \D\left(\Delta_{\nu \circ F_{\mu}^{-1}, \Lambda}\right)\subseteq C_{\nu \circ F_{\mu}^{-1},\Lambda}$ and for all $ g\in \D(A_X)$ we have $2A_X(g)=\Delta_{\nu \circ F_{\mu}^{-1}, \Lambda}(g)$ and $\nabla_{\mu} g(0)=\nabla_{\mu} g(1)=0$. Now,  Theorem \ref{thm:laplace_backward} together with Lemma \ref{lem:Equi_cont_function} yields for $f\in \D(A_Y)$
\begin{align*}
2A_Y ( f)&=2A_Y ( f)\circ \check{F}_{\mu}^{-1}\circ {F}_{\mu}=2A_X(f \circ \check{F}_{\mu}^{-1})\circ {F}_{\mu}\\
&=\Delta_{\nu \circ F_{\mu}^{-1}, \Lambda}(f\circ \check{F}_{\mu}^{-1})\circ {F}_{\mu}
=\Delta_{\nu,\mu}(f),
\end{align*}
and $\nabla_{\mu} f \circ \check{F}_{\mu}^{-1}(0)=\nabla_{\mu} f \circ \check{F}_{\mu}^{-1}(1)=0$.
Hence, the Neumann boundary conditions, given in \eqref{eq:EWC}, are 
satisfied.

\begin{Exa}
In the following, we consider the case $\mu=\nu$. The occupation formula of the local time yields for $t$ a non-negative real number
	\begin{align*}
	\Phi_t=\int L^{t}_x \d \nu \circ F_{\mu}^{-1}(x)=
	\int L^{t}_x \d \Lambda(x)
	=\int_{[0,t]} \1_{[0,1]}(B_s) \d \Lambda(s).
	\end{align*}
A simulation of a $\mu$-$\mu$-gap diffusion Brownian path when $\mu$ is the Bernoulli measure with probability weight vector $(1/2, 1/2)$ supported on the middle third Cantor set is given in Figure \ref{fig}. In order to generate a $\mu$-$\mu$-gap diffusion path we begin with a simulation of a standard Brownian path $B_{t}$, as depicted in Figure  \ref{fig}\,(\subref{a}), and the corresponding time-change function $\Phi_t$ with respect to $\Lambda$, see Figure  \ref{fig}\,(\subref{b}). From this, we can generate the associated gap diffusion path $B_{\widehat{\Phi}^{-1}_{t}}$ with speed measure $\Lambda$, as shown in Figure  \ref{fig}\,(\subref{c}). The image of $B_{\widehat{\Phi}^{-1}_{t}}$ under $\check{F}_{\mu}^{-1}$ establishes a realisation of a $\mu$-$\mu$-gap diffusion path, see Figure  \ref{fig}\,(\subref{d}). Note, despite the (Euclidean) jumps, this path is continuous with respect to the final topology $\tilde{\tau}$ on $\tilde{E}$. 
	\begin{figure}[ht!]
	\begin{subfigure}[t]{0.475\textwidth}
	\scalebox{0.85}{\input{tik}}
	\vspace*{-2em}
	\subcaption{Simulation of a Brownian path $B_{t}$.\hspace*{\fill}\label{a}}
	\end{subfigure}
	\hspace*{\fill}
	\begin{subfigure}[t]{0.475\textwidth}
	\scalebox{0.85}{\input{tik1}}
	\vspace*{-2em}
	\subcaption{The time-change function $\Phi_{t}$ with respect to the Brownian path $B_{t}$ and speed measure $\Lambda$.\label{b}}
	\end{subfigure}\\
	\begin{subfigure}[t]{0.475\textwidth}
	\scalebox{0.85}{\input{tik3}}
	\vspace*{-2em}
	\subcaption{$\Lambda$-gap diffusion path $B_{\widehat{\Phi}^{-1}_{t}}$.\hspace*{\fill}\label{c}}
	 \end{subfigure}
	\hspace*{\fill}
	\begin{subfigure}[t]{0.475\textwidth}
	\scalebox{0.85}{\input{tik5}}
	\vspace*{-2em}
	\subcaption{$\mu$-$\mu$-gap diffusion path $\check{F}^{-1}_{\mu}(B_{\widehat{\Phi}^{-1}_{t}})$.\hspace*{\fill}\label{d}}
	\end{subfigure}
	\vspace*{-0.5em}
	\caption{Simulation of a path of a generalised gap diffusion.\label{fig}}
	\end{figure}
\end{Exa}

In conclusion, 
we observe one can interpret the strong generalised Kre\u{\i}n--Felller operator $\Delta_{\nu,\mu}$ with twice the infinitesimal generator of the $\nu$-$\mu$-gap diffusion, 
and that this generalised diffusion is connected to the classical gap diffusion via the transformation given by $\check{F}^{-1}_\mu$.

\section*{Acknowledgements}

This research was supported by the DFG grant Ke 1440/3-1. The first and third authors would like to thank Institut Mittag-Leffler for their hospitality during the workshop \textsl{Thermodynamic Formalism --- Applications to Geometry, Number Theory, and Stochastics}, where part of this work was carried out.

\printbibliography

\end{document}